\documentclass[12pt]{amsart}
\usepackage{amssymb}
\begin{document}
\oddsidemargin=-1cm
\evensidemargin=-1cm
%\hoffset 5.5truemm
%\setlength{\textwidth}{25cm}

\title[Analyticity of CR mappings]{On the analyticity of CR mappings between nonminimal
hypersurfaces}
\author[Peter Ebenfelt]{Peter Ebenfelt}
\address{Department of Mathematics, 0112, University
of California at San Diego, La Jolla, CA 92093-0112} \email{
ebenfelt@math.kth.se}
%\keywords \endkeywords
\date{\number\year-\number\month-\number\day}
%\enddate
%\loadeufm

%\def\Label#1{\label{#1}{\bf \hbox{  }(#1)\hbox{  } }}
\def\Label#1{\label{#1}}
\def\1#1{\ov{#1}}
\def\2#1{\widetilde{#1}}
\def\3#1{\mathcal{#1}}
\def\4#1{\widehat{#1}}

\def\s{s}
\def\k{\kappa}
\def\ov{\overline}
\def\span{\text{\rm span}}
\def\tr{\text{\rm tr}}
\def\xo {{x_0}}
\def\Rk{\text{\rm Rk\,}}
\def\sg{\sigma}
\def \emxy{E_{(M,M')}(X,Y)}
\def \semxy{\scrE_{(M,M')}(X,Y)}
\def \jkxy {J^k(X,Y)}
\def \gkxy {G^k(X,Y)}
\def \exy {E(X,Y)}
\def \sexy{\scrE(X,Y)}
\def \hn {holomorphically nondegenerate}
\def\hyp{hypersurface}
\def\prt#1{{\partial \over\partial #1}}
\def\det{{\text{\rm det}}}
\def\wob{{w\over B(z)}}
\def\co{\chi_1}
\def\po{p_0}
\def\fb {\bar f}
\def\gb {\bar g}
\def\Fb {\ov F}
\def\Gb {\ov G}
\def\Hb {\ov H}
\def\zb {\bar z}
\def\wb {\bar w}
\def \qb {\bar Q}
\def \t {\tau}
\def\z{\chi}
\def\w{\tau}
\def\Z{\zeta}

\def \T {\theta}
\def \Th {\Theta}
\def \L {\Lambda}
\def\b {\beta}
\def\a {\alpha}
\def\o {\omega}
\def\l {\lambda}

\def \im{\text{\rm Im }}
\def \re{\text{\rm Re }}
\def \Char{\text{\rm Char }}
\def \supp{\text{\rm supp }}
\def \codim{\text{\rm codim }}
\def \Ht{\text{\rm ht }}
\def \Dt{\text{\rm dt }}
\def \hO{\widehat{\mathcal O}}
\def \cl{\text{\rm cl }}
\def \bR{\mathbb R}
\def \bC{\mathbb C}
\def \C{\mathbb C}
\def \bL{\mathbb L}
\def \bZ{\mathbb Z}
\def \bN{\mathbb N}
\def \scrF{\mathcal F}
\def \scrK{\mathcal K}
\def \scrM{\mathcal M}
\def \cR{\mathcal R}
\def \scrJ{\mathcal J}
\def \scrA{\mathcal A}
\def \scrO{\mathcal O}
\def \scrV{\mathcal V}
\def \scrL{\mathcal L}
\def \scrE{\mathcal E}
\def \hol{\text{\rm hol}}
\def \aut{\text{\rm aut}}
\def \Aut{\text{\rm Aut}}
\def \J{\text{\rm Jac}}
\def\jet#1#2{J^{#1}_{#2}}
\def\gp#1{G^{#1}}
\def\gpo{\gp {2k_0}_0}
\def\emmp {\scrF(M,p;M',p')}
\def\rk{\text{\rm rk}}
\def\Orb{\text{\rm Orb\,}}
\def\Exp{\text{\rm Exp\,}}
\def\Span{\text{\rm span\,}}
\def\d{\partial}
\def\D{\3J}
\def\pr{{\rm pr}}
\def\dbl{[\![}
\def\dbr{]\!]}

\def \D{\text{\rm Der}\,}
\def \Rk{\text{\rm Rk}\,}
\def \ima{\text{\rm im}\,}

\abstract Let $M\subset \bC^2$ be a connected real-analytic hypersurface
containing a connected complex hypersurface $E\subset \bC^2$, and let
$f\colon M\to \bC^2$ be a smooth CR mapping sending $M$ into another real-analytic
hypersurface $M'\subset \bC^2$. In this paper, we prove that if $f$
does not collapse $E$ to a point and does not collapse $M$ into the
image of $E$, and if the Levi form of $M$ vanishes to first order 
along $E$, then $f$ is real-analytic in a neighborhood of $E$. In
general, the corresponding statement is false if the Levi form of $M$
vanishes to second order or higher, in view of an example due to the
author. We also show analogous results in higher dimensions provided
that the target $M'$ satisfies a certain nondegeneracy condition. 

The main ingredient in the proof, which seems to be of independent
interest, is the prolongation of the system defining a CR mapping
sending $M$ into $M'$ to a Pfaffian system on $M$ with singularities
along $E$. The nature of the singularity is described by the order of
vanishing of the Levi form along $E$. 
\endabstract

\newtheorem{Thm}{Theorem}[section]
\newtheorem{Def}[Thm]{Definition}
\newtheorem{Cor}[Thm]{Corollary}
\newtheorem{Pro}[Thm]{Proposition}
\newtheorem{Lem}[Thm]{Lemma}
\newtheorem{Rem}[Thm]{Remark}
\newtheorem{Ex}[Thm]{Example}

\maketitle

\section{Introduction}

Let $M$ be a real-analytic hypersurface in $\bC^{n+1}$. As a
submanifold of complex space, $M$ inherits a partial complex
structure---a {\it CR structure}---in the following way (which, of
course, is well known;  see 
e.g.\ the books \cite{Bo} and \cite{BER} for basic 
concepts and notions in CR geometry).
The {\it CR bundle} $\mathcal
V$ of $M$, yielding the structure, is a rank $n$ subbundle of $\bC T
M$, the complexified tangent 
bundle on $M$, defined by $$\mathcal
V:=T^{0,1}\bC^{n+1}\cap \bC TM,$$ 
where $T^{0,1} \bC^{n+1}$ denotes the bundle of $(0,1)$ tangent
vectors on $\bC^{n+1}$. The two basic properties of $\mathcal V$ are
the following: (a) The commutator of two sections of $\mathcal V$ is
again a section of $\mathcal V$ (formal integrability); (b) $\mathcal
V\cap \bar \mathcal V$ consists of only the zero section. Sections of
the CR bundle will be called {\it CR vector fields}.

Let $L_{\bar 1},\ldots, L_{\bar n}$ denote a
basis for the CR vector fields on $M$ near a point $p_0\in M$. A
smooth mapping $f\colon M\to \bC^{k}$, with $f=(f_1,\ldots f_k)$,
is called {\it CR} (near $p_0$) if
\begin{equation}\Label{def-CRmap} L_{\bar A} f_j=0,\quad
1\leq A\leq n,\ 1\leq j\leq k.\end{equation} If $f$ sends $M$ into
another real-analytic hypersurface $M'\subset \bC^k$, i.e.\ if $f$
satisfies a nonlinear equation
\begin{equation}\Label{nonlinear} \rho'(f,\bar f)=0,\end{equation}
for some real-analytic function $\rho'(Z,\bar Z)$ with $d\rho'\neq
0$, then $f$ is CR if and only if the tangent mapping $f_*\colon
\bC T M \to \bC T M'$ sends the CR bundle $\mathcal V$ of $M$ into
the CR bundle $\mathcal V'$ of $M'$. In this paper, we shall address
the problem 
of describing conditions which imply that a CR mapping $f$ sending $M$
into $M'$ is necessarily real-analytic near $p_0$. We should point
out that a CR mapping $f$ is real-analytic near $p_0$ if and only
if $f$ extends holomorphically to a neighborhood of $p_0$ in
$\bC^{n+1}$ (see e.g.\ \cite{BER}, Chapter I).

There is an extensive literature on this subject if $M$ is assumed
to be {\it minimal} at $p_0$, i.e. if there are no complex
hypersurfaces through $p_0$ in $\bC^{n+1}$ contained in $M$. In
this case, all CR mappings extend holomorpically to one side of
$M$ (see \cite{Tr}; cf.\ also \cite{BT}, \cite{Tu}) and results
on the problem decribed above are often referred to as reflection
principles. We mention here the papers \cite{Pinch}, \cite{Lewy}, \cite{BJT},
\cite{BR1}, \cite{DF}, \cite{Me}, \cite{Hu1},  and
refer the reader to the notes in \cite{BER}, Chapter IX, and the
survey article \cite{Hu2} for a more detailed history. We should
point out, however, that most results mentioned above require
additional nondegeneracy conditions on the manifolds and mappings
which do not appear to be necessary conditions. Thus, even though much is
known, the problem is still not completely resolved in the minimal
case.

In the present paper, we shall assume that $M$ is not minimal at
$p_0$; i.e. there exists a complex hypersurface
$E\subset\bC^{n+1}$ (which we shall assume is connected) through
$p_0$ contained in $M$. To state our main results, we need to
define an invariant which measures the order of vanishing of the
Levi form along $E$. To be more precise, we must introduce some
notation. It is well known (\cite{CM}; see also \cite{BER},
Chapter IV]) that, for any real-analytic hypersurface
$M\subset\bC^{n+1}$ and $p_0\in M$, there are local holomorphic
coordinates $(z,w)\in \bC^n\times\bC$ vanishing at $p_0$ such that
$M$ is defined locally near $p_0=(0,0)$ by an equation of the form
\begin {equation}\Label{def-eq}\im w=\phi(z,\bar z,\re
w),\end{equation} where $\phi(z,\chi,s)$ is a holomorphic function
satisfying
\begin
{equation}\Label{def-normal}\phi(z,0,s)=\phi(0,\chi,s)=0.\end{equation}
Such coordinates are called {\it normal coordinates} for $M$ at
$p_0$. Although they are not unique, one can define an integer
$m\geq 0$, which is independent of the choice of normal
coordinates, as follows. If $\phi(z,\bar z,s)\equiv 0$, we set
$m=\infty$; otherwise, we define $m$ by Taylor expanding in
$s$,\begin {equation}\Label {def-m}
\phi(z,\chi,s)=\phi_m(z,\chi)s^m+O(s^{m+1}),\end{equation} such
that $\phi_m(z,\chi)\not\equiv 0$.  The fact that the integer $m$
is a biholomorphic invariant of $M$ at $p_0$ was first observed in
\cite{Me}. If $m=\infty$, i.e. if $\phi(z,\bar z, s)\equiv 0$,
then $M$ is said to be {\it Levi flat}.

When $M$ is not minimal at $p_0$, the complex hypersurface
$E\subset M$ is given, in normal coordinates, by $w=0$ and the
fact that $M$ is not minimal at $p_0$ is therefore characterized
by $m\geq 1$. It is not difficult to see that the integer $m$,
which {\it a priori} depends on $p_0\in E\subset M$, is constant
along $E$. Indeed, as mentioned above, $m$ measures the order of
vanishing of the Levi form of $M$ along the complex hypersurface
$E$ (see Proposition \ref{m-hAB} for the precise statement). We
shall say that $M$ is of {\it $m$-infinite type} along $E$ (or at
$p_0$); recall that for a real-analytic hypersurface $M$, being
nonminimal at $p_0$ is equivalent to being of infinite type in the
sense of Kohn \cite{Kohn} and Bloom--Graham \cite{BG} at $p_0$.  One
of our main results,  for 
$n=1$, is the following.

\begin{Thm}\Label{Main-1} Let $M\subset \bC^{2}$ be a real-analytic, connected
hypersurface which is of $1$-infinite type along a connected
complex curve (or, equivalently, hypersurface) $E\subset \bC^2$
contained in $M$. Let $f\colon M\to \bC^{2}$ be a
$C^\infty$-smooth CR mapping. Assume that $f(M)$ is contained in a
real-analytic hypersurface $M'\subset \bC^{2}$, and that
\begin{enumerate}
\item[(i)] $f(M)\not\subset f(E)$;
\item[(ii)] $f|_E$ is not constant.
\end{enumerate}
Then, $f$ extends as a holomorphic mapping from an open
neighborhood of $E$ in $\bC^2$, i.e.\ there exists a holomorphic
mapping $F\colon U\to \bC^{2}$, where $U$ is an open neighborhood
of $E$ in $\bC^2$, such that $F|_M=f$.
\end{Thm}

Theorem \ref{Main-1} will follow from a more general, but slightly
more technical result which is valid in any number of dimensions.
However, before explaining this more general result we make a few
remarks. First, we would like to point out that Theorem
\ref{Main-1} fails in general if $M$ is of $m$-infinite type along
$E$ with $m\geq 2$ as the following example from \cite{E-JGA}
shows.

\begin{Ex}\Label{Ex-1} {\rm Let $M\subset \bC^2$ be defined by
\begin{equation}\Label{Ex12-1}
\im w=\theta(\arctan |z|^2,\re w),\end{equation} where
$t=\theta(\xi,s)$ is the unique solution of $\xi(s^2+t^2)-t=0$
with $\theta(0,0)=0$. One can show that $M$ is real-analytic near
$(z,w)=(0,0)$ and of $2$-infinite type along the complex curve
$E:=\{w=0\}$ which is contained in $M$. The restriction $f$ to $M$
of the mapping $F(z,w):=(z,g(w))$, where $g(w)=e^{-1/w}$ if $\re
w>0$ and $g(w)=0$ if $\re w\leq 0$, is a $C^\infty$-smooth CR
mapping which does not extend holomorphically to any neighborhood
of $(0,0)$. (Indeed, the second component of $f$ vanishes to
infinite order along $E$.) Moreover, it is shown in \cite{E-JGA}
that $f(M)$ is contained in the real-analytic (indeed,
real-algebraic) hypersurface $M'\subset \bC^2$ defined by
\begin{equation}\Label{Ex12-2} \im w=(\re w)|z|^2.\end{equation}
Observe that conditions (i) and (ii) of Theorem \ref{Main-1} are
both satisfied.} \end{Ex}

The CR mapping in Example \ref{Ex-1} does not extend to either
side of the hypersurface $M$ near $E$. Indeed, the
conclusion of Theorem \ref{Main-1}, with $M$ of $m$-infinite type
for any $1\leq m<\infty$ and $f$ merely assumed to be continuous,
was proved in \cite{EH} under the {\it additional assumption} that
the CR mapping extends holomorphically to one side of $M$. As
shown by Example \ref{Ex-1}, the conclusion fails for $m\geq 2$
without assuming one-sided extension. It is also noteworthy that,
without this assumption, the conclusion fails with
$C^\infty$-regularity replaced by $C^k$-regularity for any finite
$k$ (in contrast, again, with the result in \cite{EH}) as is shown
by the following example from \cite{EH}.

\begin{Ex}\Label{Ex-2} {\rm Let $M\subset \bC^2$ be defined by
\begin{equation}\Label{Ex2-1}
 \im w=(\re w)|z|^2,\end{equation} and $M'_k\subset\bC^2$ by
\begin{equation}\Label{Ex2-2} \im w=(\re w)h_k(z,\bar{z}),\end{equation}
where $h_k$ is defined as follows. Let
\begin{equation}\Label{Ex2-3}\begin{aligned}
\phi_k(s,t)&\, :=\re((s+it)^k)=s^k+\sum_{j=0}^{k-1}
a_js^jt^{k-j}\\ \psi_k(s,t)&\, :=\im((s+it)^k)=\sum_{j=0}^{k-1}
b_js^jt^{k-j}\end{aligned}\end{equation} and define
\begin{equation}\Label{Ex2-4}h_k(z,\bar z):=\frac{\sum_{j=0}^{k-1}
b_j|z|^{2(k-j)}}{1+\sum_{j=0}^{k-1}
a_j|z|^{2(k-j)}},\end{equation} where the $a_j$ and $b_j$ are
defined by (\ref{Ex2-3}). The mapping $F(z,w):=(z,g(w))$, with
$g(w)=-w^k$ for $\re w\le 0$ and $g(w)=w^k$ for $\re w>0$ is a
$C^{k-1}$-differentiable CR mapping from $M$ into $M'_k$, which
satisfies (i) and (ii) of Theorem \ref{Main-1} but does not extend
holomorphically to a neighborhood of $0$ in $\bC^2$. Observe that $M$
is of $1$-infinite type along $w=0$.}\end{Ex}

In order to state the more general result for real-analytic
hypersurfaces in higher dimensional complex space, we need the
following notion from \cite{Me}. Assume that $M$ is not Levi flat
and write the defining equation (\ref{def-eq}) of $M$ near
$p_0=(0,0)$ as
\begin{equation}\Label{def-eq2}\im w=(\re w)^m\psi(z,\bar z,\re w),
\end{equation}
where $m$ is the integer defined by (\ref{def-m}), and also write
\begin{equation}\Label{psi-exp} \psi(z,\chi,0)=\sum_{\alpha}
a_\alpha(z)\chi^\alpha.\end{equation} Then, we shall say that $M$
is {\it $m$-essential} at $p_0=(0,0)$ if the ideal generated by
the collection $a_\alpha(z)$, $\alpha\in \mathbb Z_+^n$, is of
finite codimension in the ring $\bC\{z\}$ of convergent power
series. This notion is analogous to that of essential finiteness in
the finite type case and 
was shown, in \cite{Me}, to be independent of
the choice of normal coordinates for $M$ at $p_0$. If $M$ is of
$m$-infinite type at $p_0$, for some $1\leq m<\infty$, and $M$ is
$m$-essential at $p_0$, then we shall also say that $M$ is {\it
weakly essential} at $p_0$.

\begin{Thm}\Label{Main-2} Let $M\subset \bC^{n+1}$ be a real-analytic,
connected hypersurface which is of $1$-infinite type along a
connected complex  hypersurface $E\subset \bC^{n+1}$ contained in
$M$. Let $f\colon M\to \bC^{n+1}$ be a $C^\infty$-smooth CR
mapping. Assume that $f(M)$ is contained in a real-analytic
hypersurface $M'\subset \bC^{n+1}$, and that
\begin{enumerate}
\item[(i)] $f(M)\not\subset f(E)$;
\item[(ii)] $f|_E\colon E\to \bC^{n+1}$ is a finite mapping.
\end{enumerate}
Assume, in addition, that $M'$ is weakly essential at some point
of $f(E)\subset M'$. Then, $f$ extends as a holomorphic mapping
from an open neighborhood of $E$ in $\bC^{n+1}$, i.e.\ there
exists a holomorphic mapping $F\colon U\to \bC^{n+1}$, where $U$
is an open neighborhood of $E$ in $\bC^{n+1}$, such that $F|_M=f$.
\end{Thm}

The conclusion of Theorem \ref{Main-2}, under the additional
assumption that $f$ extends holomorphically to one side of $M$,
was proved in \cite{HMM} for $M$ of $m$-infinite type, with $1\leq
m<\infty$, and $f$ merely assumed to be $C^1$-smooth.

We would also like to mention some results, related to those in this paper,
without describing them in detail. The first is the
Baouendi-Rothschild reflection principle \cite{BR2} which deals
with $C^\infty$-smooth CR mappings which are assumed to extend to one
side of a Levi-nonflat real-analytic
hypersurface in $\bC^2$, and the second is a general result
(\cite{BHR}) for $C^\infty$-smooth CR mappings of real algebraic,
holomorphically nondegenerate submanifolds in any dimension. We
also mention that there are some results (see e.g.\ \cite{Forst}, \cite{Han2},
\cite{Hay},  \cite{Lamel}) for CR mappings between manifolds in 
spaces of different dimension.

The idea of the proof of our principal result, Theorem \ref{Main-2}, is roughly
the following. We first prolong the system (\ref{def-CRmap}),
(\ref{nonlinear}) (or, more precisely, its intrinsic counterpart) to a
Pfaffian system on $M$. The latter system develops a singularity along
the complex hypersurface $E\subset M$, and the nature of this
singularity near points in
general position on $E\subset M$ is described by the invariant $m$
(see Theorem \ref{thm-maintech}); when $m=1$, the 
singularity is regular in a certain sense (Fuchsian), and one can use
known results about such systems combined with the Hanges--Treves
propagation theorem to complete the proof of Theorem 1.4. The details
are carried out below.

\section{Singular Pfaffian systems for CR mappings}\Label{MTR}

As above, let $M$ be a real-analytic hypersurface  in $\bC^{n+1}$
and $p_0\in M$. We shall keep the notation introduced in the
previous section. Thus, we assume that $(z,w)\in \bC^n\times \bC$ are
normal coordinates for $M$ at $p_0=(0,0)$ and that $M$ is defined
by an equation of the form (\ref{def-eq}), where $\phi(z,\chi,s)$
is a holomorphic function satisfying (\ref{def-normal}). We also
assume that $M$ is of $m$-infinite type, for some $m\geq 1$, along
the complex hypersurface $E$, where $m$ is defined by (\ref
{def-m}). It is also shown in \cite{Me} that the lowest order $r$
in the Taylor expansion
\begin{equation}\Label{def-r}
\phi_m(z,\chi)=\sum_{|\alpha|+|\beta|^\geq r}c_{\alpha\beta}
z^\alpha\chi\beta,\end{equation} where $\phi_m(z,\chi)$ is defined
by (\ref{def-m}), is a biholomorphic invariant of $M$ at $p_0$.
With these definitions of the integers $m$ and $r$, we shall say
that $M$ is of {\it $m$-infinite type $r$} at $p_0=(0,0)$.

Let us write the defining equation of $M$ as in (\ref{def-eq2}).
We shall also say that $M$ is {\it $m$-infinite
$\ell$-nondegenerate} at $p_0=(0,0)$ if
\begin{equation}\Label{def-infnondeg}\span_\bC\left\{\frac{
\partial^{|\alpha|}}
{\partial \bar z^\alpha}\left(\frac{\partial\psi}{\partial
z}\right)\colon\forall |\alpha|\leq
\ell\right\}=\bC^n,\end{equation} where $\partial/\partial
z=(\partial/\partial z_1,\ldots,\partial/\partial  z_n)$ and
standard multi-index notation is used. The notion is completely
analogous to that of $\ell$-nondegeneracy for hypersurfaces of
finite type (see e.g.\ \cite{BER}, Chapter XI), and the reader can
verify, as in the finite type case, that the definition above is
independent of the choice of normal coordinates. We should point
out, as is remarked in \cite{HMM}, that if $M$ is $m$-essential at
some point on $E\subset M$, then $M$ is $m$-infinite
$n$-nondegenerate outside a proper complex subvariety of $E$.

Given two smooth manifolds $M$, $M'$ of the same dimension $2n+1$,
let us denote by $J^k(M,M')_{(p,p')}$ the space of $k$-jets at
$p\in M$ of smooth mappings $f\:M\to M'$ with $f(p)=p'\in M'$.
Given coordinate systems $x=(x_1,\ldots, x_{2n+1})$ and
$x'=(x'_1,\ldots, x'_{2n+1})$ on $M$ and $M'$ near $p$ and $p'$,
respectively, there are natural coordinates
$\lambda^k:=(\lambda^\beta_i)$, where $1\leq i\leq 2n+1$ and
$\beta\in\mathbb Z_+^{2n+1}$ with $1\leq |\beta|\leq k$, on
$J^k(M,M')_{(p,p')}$ in which the $k$-jet at $p$ of a smooth
mapping $f\:M\to M'$ is given by
$\lambda_i^\beta=(\partial_x^\beta f_i)(p)$, $1\leq |\beta|\leq k$
and $1\leq i\leq 2n+1$. The main technical result in this paper is
the following.

\begin{Thm}\Label{thm-maintech} Let $M,M'\subset \bC^{n+1}$ be
real-analytic hypersurfaces which are of $m$-infinite and
$m'$-infinite type respectively, for some integers $m,m'\geq 1$,
along complex hypersurfaces $E\subset M$ and $E'\subset M'$. Let
$p_0\in E\subset M$ and $p_0'\in E'\subset M'$. Assume that $M'$
is $m'$-infinite $\ell$-nondegenerate at $p'_0\in M$ and $M$ is of
$m$-infinite type $2$ at $p_0\in M$. Let $f^0\colon M\to M'\subset
\bC^{n+1} $ be a $C^\infty$-smooth CR mapping such that
$f^0(p_0)=p_0'$ and such that $f=f^0$ satisfies:
\begin{enumerate}
\item[(i)] $f(M)\not\subset E'$;
\item[(ii)] $f|_E\colon E\to \bC^{n+1}$ is a local immersion at $p_0$.
\end{enumerate}
Choose local coordinates $y=(x,s)\in\bR^{2n}\times \bR$ on $M$
near $p_0$ such that $E$ is given by $s=0$. Then, for any
multi-index $(\alpha,p)\in\mathbb Z_+^{2n}\times\mathbb Z$ with
\begin{equation}\Label{maintech-0}|\alpha|+p=2\ell+2
\end{equation} and any $j=1,\ldots,2n+1$,
there is a real-analytic function $r^{\alpha,p}_j(\lambda^{k};y')(y)$
on $U$, where $k:= 2\ell+1$ and $U\subset
J^{k}(M,M')_{(p_0,p_0')}\times M\times M'$ is an open neighborhood
of $(((s^m\partial_s)^q\partial_x^\beta f^0)(p_0),f^0(p_0),p_0)$,
such that \begin{equation}\Label{maintech-1}
(s^m\partial_s)^p\partial_x^\alpha
f_j=r^{\alpha,p}_j((s^m\partial_s)^q\partial_x^\beta
f;f),\end{equation} where $1\leq |\beta|+q\leq k$, for every
smooth CR mapping $f\colon V\to M'$, where $V\subset M$ is some
open neighborhood of $p_0$, which satisfies {\rm (i)--(ii)} above
and with $(((s^m\partial_s)^q\partial_x^\beta
f)(p_0),f(p_0),p_0)\in U$; here, $y'=(y'_1,\ldots, y'_{2n+1})$ is
any local coordinate system on $M'$ near $p_0'=f^0(p_0)$ and
$f_i:=f\circ y'_i$. The functions $r^{\alpha,p}_j$ depend only on $M$,
$M'$ and the $k$-jet of $f^0$ at $p_0$. In addition, the functions
$r^{\alpha,p}_j$ are rational in $\lambda^{k}\in
J^{k}(M,M')_{(p_0,p_0')}$.
\end{Thm}

Before turning to the proof of Theorem \ref{thm-maintech}, we
mention that a similar result was proved at minimal points (i.e.\
$m=0$) in \cite{Hay}, using ideas of \cite{Han1}, \cite{Han2}, and
in \cite{E-UCP}, using a more intrinsic approach which allowed
its extension to the case of merely smooth hypersurfaces. In
this paper, we shall follow the latter approach. As a consequence,
Theorem \ref{thm-maintech} also follows with $C^\infty$-smoothness
(for the hypersurfaces $M$ and $M'$ as well as for the functions
$r^{\alpha,p}_j$) replacing real-analyticity. However, the author
has not found any direct application of this result, even though
it is clear that it does restrict the behaviour of
the mappings as $s\to 0$.

\section{Preliminaries}\Label{prelim}

The proof of Theorem \ref{thm-maintech} is similar to that of
Theorem 2 in \cite{E-UCP}. Indeed, the proof of Theorem
\ref{thm-maintech} will essentially reduce to that in \cite{E-UCP}
by suitably modifying and adapting the setup in \cite{E-UCP} to
the present situation where the hypersurfaces are of infinite
type. We first need to relate the notions of $m$-infinite type and
$m$-infinite $\ell$-nondegeneracy defined in the previous section
to the CR geometry of $M$.

Recall that $\mathcal V:=\bC TM\cap T^{0,1}\bC^{n+1}$ denotes the
CR bundle on $M$. We shall denote by $T^0 M$ the characteristic
bundle $(\mathcal V\oplus\bar\mathcal V)^\perp\subset \bC T^* M$,
and by $T'M$ the holomorphic bundle $\mathcal V^\perp\subset \bC
T^* M$. Real non-vanishing sections of $T^0M$ are called
characteristic forms, sections of $T'M$ holomorphic forms, and
sections of $\mathcal V$ CR vector fields. The reader is referred
to \cite{BER} for basic definitions and facts regarding CR
manifolds and structures.

The Levi form of $M$ at $p\in M$ is a multi-linear mapping
$\Lambda_p\colon \mathcal V_p\times\bar\mathcal V_p\times
T^0_pM\to \bC$ (or, equivalently, a tensor in $\mathcal
V_p^*\times\bar\mathcal V_p^*\times (T^0_pM)^*$) defined by
\begin {equation}\Label {def-levi}
\Lambda_p(X_p,Y_p,\theta_p):=\frac{1}{2i}\left<\theta,[X,Y]\right>_p=
-\frac{1}{2i}\left<d\theta,X\wedge Y\right>_p,\end{equation} where
$X\in\Gamma(M,\mathcal V)$ and $Y\in\Gamma(M,\bar \mathcal V)$ are
vector fields extending $X_p$ and $Y_p$, respectively, and
$\theta$ a characteristic form extending $\theta_p$. We shall
assume that $M$ is of $m$-infinite type along a complex
hypersurface $E\subset M$. We first claim that $m$ is also the
order of vanishing of the Levi form $\Lambda$ (as a function of
$p\in M$ with values in $\mathcal V_p^*\times\bar\mathcal
V_p^*\times (T^0_pM)^*$) along $E$. More precisely, we choose,
near a given point $p_0\in E\subset M$, a basis for the CR vector
fields $L_{\bar 1},\ldots, L_{\bar n}$, set
$L_{A}:=\overline{L_{\bar A}}$, choose a nonvanishing
characteristic form $\theta$, and represent $\Lambda$ by the
$n\times n$ matrix $(h_{\bar A B})$, $1\leq A,B\leq n$, where
\begin{equation}\Label{def-hAB} h_{\bar A B}(p):=-2i\Lambda_p(L_{\bar
A},L_B,\theta).\end{equation} Then the following holds.
\begin{Pro}\Label{m-hAB} If $M\subset \bC^{n+1}$ is a real-analytic hypersurface
which is of $m$-infinite type along a complex hypersurface
$E\subset M$ through $p_0\in M$, then there exists a
real-analytic, $n\times n$ matrix valued function $(h^0_{\bar
AB})$, $1\leq A,B\leq n$, in a neighborhood of $p_0$ such that the
restriction $(h^0_{\bar AB})|_E$ is not identically 0, and
$$h_{\bar AB}=\delta^mh^0_{\bar AB},$$ where $\delta$ denotes the
distance (in the Riemannian metric on $M$ inherited from the
ambient space) of $p$ to $E$ and $h_{\bar AB}$ is an $n\times n$
matrix representing the Levi form as explained above. In addition,
$M$ is of $m$-infinite type $2$ at $p_0$ if and only if $h^0_{\bar
AB}(p_0)\neq 0$ for some $1\leq A,B\leq n$.
\end{Pro}
Before proving Proposition \ref{m-hAB}, we shall introduce a
special choice of basis for the CR vector fields on $M$ near
$p_0$. Let $(z,w)$ be normal coordinates for $M$ at $p_0$, so that
$M$ is defined near $p_0=(0,0)$ by (\ref{def-eq}). We may then
take $(z,s)\in \bC^n\times \bR$, with $s=\re w$, to be local
coordinates on $M$ near $(0,0)$, and choose
\begin{equation}\Label{CRvf} L_{\bar A}=\partial_{\bar z_A}-
\frac{i\phi_{\bar z_A}}{1+i\phi_s}\partial_s,\end{equation}
 where $\phi=\phi(z,\bar z,s)$, and where we have used the
notation $\partial_{\bar Z_A}=\partial/\partial \bar z_A$ and
$\phi_s=\partial \phi/\partial s$, etc.\ Observe for future
reference that the vector fields $L_A$ and $L_B$, for any $1\leq
A,B\leq n$, commute. In the coordinates $(z,s)$, the distance
$\delta$ is comparable to $s$. We also denote by $T$ the vector
field $\partial_s$, so that $T,L_1,\ldots,L_n,L_{\bar 1},\ldots,
L_{\bar n}$ spans $\bC T M$ near $p_0$. In addition, we choose the
characteristic form $\theta$ so that $\left<\theta,T\right>=1$.
\begin{proof}[Proof of Proposition $\ref{m-hAB}$] It is not difficult to
see that the order of vanishing of the Levi form along $E$ is
independent of the choice of basis $L_{\bar 1},\ldots, L_{\bar n}$
of the CR vector fields, and characteristic form $\theta$ near
$p_0$. Thus, it suffices to prove Proposition \ref{m-hAB} using
the special choices introduced above. Since $M$ is assumed to be
of $m$-infinite type along $E$, we may write $$\phi(z,\bar
z,s)=s^m(\alpha(z,\bar z)+O(|z|^{r+1}))+O(s^{m+1}),$$ where
$\alpha(z,\bar z)\not\equiv 0$ is a homogeneous polynomial of some
degree $r\geq 2$ with $\alpha(z,0)\equiv \alpha(0,\bar z)\equiv 0$. In
particular, the $n\times n$ matrix $(\alpha_{\bar z_Az_B})$ is not
identically 0. A straightforward calculation shows that $$h_{\bar
A B}=s^m(\alpha_{\bar z_Az_B}(z,\bar
z)+O(|z|^{r-1}))+O(s^{m+1}),$$ which completes the proof of
Proposition \ref{m-hAB}.
\end{proof}

Next, for $A_1,\ldots, A_k\in\{1,\ldots, n\}$, we define,
following \cite{E-JDG}, \cite{E-UCP}, the functions
\begin{equation}\Label{def-hgeneral} h_{\bar A_1\ldots \bar A_k
D}:=\left<\mathcal L_{\bar A_k}\ldots\mathcal L_{\bar
A_1}\theta,L_D\right>,\end{equation} where $\mathcal L_{\bar
A}\omega:=\mathcal L_{L_{\bar A}}\omega$, for a holomorphic form
$\omega$, denotes the Lie derivative $L_{\bar A}\lrcorner d\omega$
along the CR vector field $L_{\bar A}$ (which is again a
holomorphic form). It is shown in \cite{E-UCP} that
\begin{equation}\Label{Lhgeneral} h_{\bar A_1\ldots \bar A_k \bar
C D}=L_{\bar C} h_{\bar A_1\ldots \bar A_kD}+h_{\bar A_1\ldots
\bar A_k} h_{\bar C D},\end{equation} where $$ h_{\bar A_1\ldots
\bar A_k}:=\left<\mathcal L_{\bar A_k}\ldots\mathcal L_{\bar
A_1}\theta,T\right>.$$ By using Proposition \ref{m-hAB} and
(\ref{Lhgeneral}) inductively, we conclude that
\begin{equation}\Label{nothing} h_{\bar A_1\ldots \bar A_k
D}=s^mh^0_{\bar A_1\ldots \bar A_k D},\end{equation} where
$h^0_{\bar A_1\ldots \bar A_k \bar D}$ is a real-analytic function
satisfying the identity
\begin{equation}h^0_{\bar A_1\ldots \bar A_k\bar C \bar D}=
L_{\bar C} h^0_{\bar A_1\ldots \bar A_k \bar D}+a_{\bar C}
h^0_{\bar A_1\ldots \bar A_k \bar D}+h_{\bar A_1\ldots \bar A_k}
h^0_{\bar C D};\end{equation} here, $a_{\bar C}$ is the function
$m(L_{\bar C}s)/s$, which is real-analytic since $L_{\bar C} s$
vanishes on $E$.

Now, define the filtration \begin{equation}\Label{filtration}
\bar\mathcal V_0=F_0(0)\supset F_1(0)\supset\ldots\supset
F_k(0)\supset\ldots \supset\{0\},\end{equation} where each
subspace $F_k(0)$ is defined as follows,
\begin{multline}\Label{Fk} F_k(0):=\\
\left\{X_0=a^DL_D(0)\in\bar\mathcal V_0\colon a^D h^0_{\bar
A_1\ldots \bar A_j D}(0)=0,\ \forall A_1\ldots A_j\in \{1,\ldots
n\},\ j\leq k\right\}. \end{multline} In (\ref{Fk}), we have used
the summation convention, i.e.\ an index appearing both as a
super- and subscript is summed over. Moreover, in what follows,
capital Roman indices ($A,B$, etc.\  ) will run over the set
$\{1,\ldots,n\}$. As in \cite{E-JDG}, one can check that $$
(X_1,\ldots,X_k,Y,\theta)\to \frac{1}{s^m}\lim_{(z,s)\to (0,0)}
\left<\mathcal L_{X_k}\ldots\mathcal L_{X_1}\theta,Y\right>,$$
defines a multi-linear mapping $$\underset{\text{\rm $k$
times}}{\underbrace{\mathcal V_0\times\ldots\times \mathcal
V_0}}\times F_{k-1}(0)\times T^0_0M\to \bC,$$ which is symmetric
in the first $k$ positions. Set $r_k=n-\dim F_k(0)$. By a constant
linear change of the $L_A$, we may adapt the basis $L_A$ of the
complex conjugate CR vector fields to the filtration
(\ref{filtration}) so that $L_{r_k+1}(0),\ldots L_n(0)$ spans
$F_k(0)$ for each $k=0,1,\ldots, \ell$, where $\ell$ is the
smallest integer for which $F_\ell(0)$ is minimal. We also adopt
the index convention from \cite{E-UCP}. For $j=1,2, \ldots$, Greek
indices $\alpha^{(j)},\beta^{(j)}$, etc., will run over the set
$\{1,\ldots, r_{j-1}\}$ and small Roman indices $a^{(j)},
b^{(j)}$, etc.,  over $\{r_{j-1}+1,\ldots, n\}$. As mentioned
above, capital Roman indices $A,B$, etc., will run over
$\{1,\ldots, n\}$. The linear change of the $L_A$ which adapts the
basis to the filtration (\ref{filtration}) corresponds to a linear
change of the coordinates $z$ in normal coordinates, so that we
may still assume that the basis $L_{\bar A}$ is of the form
(\ref{CRvf}).

A straightforward calculation in normal coordinates (cf.\
\cite{E-SJM}) shows that $M$ is $m$-infinite $\ell$-nondegenerate
if and only if $F_{\ell}(0)=\{0\}$ and $F_j(0)\neq \{0\}$ for
$j<\ell$. This is equivalent to the fact that there are indices
$\underline A^j:=A^j_1\ldots A^j_{k_j}$, $j=1,\ldots,n$ with
$k_j\leq \ell$, such that the $n\times n$ matrix
$(h^0_{\underline{A^j} D}(0))$, $1\leq j,D\leq n$ is invertible.
In addition, by the choice of basis $L_A$ and the index
convention, we have the following two basic facts, whose proofs
are elementary and left to the reader (c.f.\ also
\cite{E-UCP}). First, we have
\begin{equation}\Label{fact-1} h^0_{\bar A_1\ldots\bar A_j a^{(k)}}
(0)=0,\quad\forall j<k,\end{equation} and, secondly, 
\begin{Lem}\Label{fact-2} If $(v^A)\in \bC^n$ satisfies $$v^{a^{(k)}} h^0_
{\bar A_1\ldots\bar A_k a^{(k)}}(0)=0,\quad \forall A_1\ldots
A_k\in \{1,\ldots,n\},$$ for some $k\leq \ell$, then
$v^{\alpha^{(k+1)}}=0$.
\end{Lem}
For future reference, we also record here the fact that any
commutator $[X,Y]$, for $X,Y\in\{T,L_A,L_{\bar A}\}$, is a
multiple of $T$, and a straightforward calculation shows that
\begin{equation}\Label{fact-3} \begin {aligned}
\, [L_{\bar A},s^m T]\,  &=
    ms^{m-1}(L_{\bar A}s) T+s^m [L_{\bar A},T]\\
\, &=-s^m h^0_{\bar A},\end{aligned}\end{equation} where
$h^0_{\bar A}$ is a real-analytic function near $0$ since $L_{\bar
A} s$ vanishes on $E$.

\section{Proof of Theorem
$\ref{thm-maintech}$}\Label{proof-maintech}

We shall keep the notation and conventions introduced in previous
sections for the real-analytic hypersurface
$M\subset\bC^{n+1}$. We shall also need the real vector field
$S=s^m T$, where $m$ is the invariant associated to $M$ in Theorem
\ref{thm-maintech}; i.e.\ $M$ is of $m$-infinite type along $E$. For
convenience of  
notation, we shall denote the target hypersurface $M'\subset
\bC^{n+1}$ by $\hat M$ and denote corresponding objects for $\hat
M$ by placing a hat over them; i.e.\ $\hat E$ denotes the complex
hypersurface through $\hat p_0$ contained in $\hat M$, $(\hat z,
\hat s)\in \bC^n\times \bR$ denote local coordinates on $\hat M$
near $\hat p_0=(0,0)$, $\hat L_{\bar A}$ denotes a basis for the
CR vector fields on $\hat M$ of the form (\ref{CRvf}), $\hat T$
denotes the real vector field $\partial/\partial\hat s$,  etc.\

Assume that $f\colon M\to \hat M$ is a smooth ($C^\infty$) CR
mapping defined near $p_0=0$ in $M$ such that $f(0)=\hat p_0=0$
and $f$ satisfies conditions (i)--(ii) in Theorem
\ref{thm-maintech}. Recall that a smooth mapping $f\colon M\to
\hat M$ is called CR if $f_*(\mathcal V_p)\subset \hat\mathcal
V_{f(p)}$, where $f_*\colon \bC TM\to \bC T\hat M$ denotes the
tangent mapping or push forward, for every $p\in M$. 
When $E, \hat E\subset \bC^{n+1}$ are complex hypersurfaces
contained in $M$ and $\hat M$, respectively, then any CR mapping
$f\colon M\to \hat M$ sends $E$ into $\hat E$. Thus, condition
(ii) in Theorem \ref{thm-maintech} is equivalent to $f|_E\colon
E\to \hat E$ being a local diffeomorphism. Also, observe that if
$f|_E\colon E\to \hat E$ is a local diffeomorphism near $0$ then $
f_*(\mathcal V_0)= \hat\mathcal V_{0}$. We introduce the smooth
$GL(\bC^n)$-valued function $(\gamma^A_B)$, smooth complex-valued
functions $\eta^A$, and real-valued function $\xi$ so that
\begin{equation}\Label{push-1}
f_*(L_B)=\gamma^A_B\hat L_A,\quad f_*(L_{\bar
B})=\overline{\gamma^A_B}\hat L_{\bar A},\quad f_*(S)=\xi\hat
S+\eta^A\hat L_A+\overline{\eta^A}\hat L_{\bar A}.\end{equation}
Observe that $\xi$ is {\it a priori} possibly singular along
$f^{-1}(\hat E)$
since $\hat S$ vanishes along $\hat E$. Indeed, using the
coordinates introduced in the previous section, we have
\begin{equation}\Label{f-xi}
s^m\frac{\partial \hat s}{\partial s}=\hat s^{\hat
m}\xi+\frac{i}{1-i\phi_s}\frac{\partial \phi}{\partial z_A}
\frac{\partial \hat z_A}{\partial s}-\frac{i}{1+i\phi_s}\overline
{\frac{\partial \phi}{\partial z_A}\frac{\partial \hat
z_A}{\partial s}}.
\end{equation}
However, we shall show (Proposition \ref{xi}) that if $M$ is of
$m$-infinite type $2$ at $0$, then $\xi$ 
is in fact smooth near $0$. 

We can write (\ref{push-1}) using matrix notation as
\begin{equation}\Label{push-2} f_*(S,L_B,L_{\bar B})=(\hat S,\hat L_A,\hat L_{\bar
A})\pmatrix \xi&0&0\\\eta^A&\gamma^A_B&0\\\overline{\eta^A}&0
&\overline{\gamma^A_B}\endpmatrix .\end{equation} Since $f$
satisfies condition (i) in Theorem \ref{thm-maintech}, it is well
known that in fact $f(M\setminus E)\subset \hat M\setminus \hat
E$. Hence, if we let $\theta,\theta^A, \theta^{\bar A}$, where
$\theta$ is as in section \ref {prelim}, be a dual basis (of
$1$-forms) to $T,L_A,L_{\bar A}$ then, by duality, we have
(outside $E$)
\begin{equation}\Label{pull-1}
f^*\pmatrix \hat\theta/\hat s^{\hat m}\\\hat
\theta^A\\\hat\theta^{\bar A}\endpmatrix =\pmatrix
\xi&0&0\\\eta^A&\gamma^A_B&0\\\overline{\eta^A}&0
&\overline{\gamma^A_B}\endpmatrix \pmatrix \theta/s^m\\ \theta^B\\
\theta^{\bar B}\endpmatrix .\end{equation} We shall make use of
the two identities \begin{equation}\Label{id-1} \left<df^*\hat
\omega,X\wedge Y\right>= \left<d\hat \omega,f_*X\wedge
f_*Y\right>,\end {equation} where the left side is evaluated at
$p\in M$ and the right side at $f(p)$,  which holds for any 1-form
$\hat \omega$ on $\hat M$ and vector fields $X$, $Y$ on $M$, and
also
\begin{equation}\Label{id-2} \left<d \omega,X\wedge Y\right>=
-\left<\hat \omega,[X,Y]\right>,\end{equation} which holds for any
1-form $\omega\in \{\theta/s^m,\theta^A, \theta^{\bar A}\} $ and
vector fields $X,Y\in\{S,L_A,L_{\bar A}\}$ on $M$ since
$\theta/s^m,\theta^A, \theta^{\bar A}$ is a dual basis (outside
$E$) to $S,L_A,L_{\bar A}$. First, we apply (\ref{id-1}) with
$\hat\omega=\hat \theta/\hat s^{\hat m}$, $X=L_{\bar A}$, and
$Y=L_B$, and obtain
\begin{equation}\begin{aligned}\Label{calc-11} \left<df^*(\hat
\theta/\hat
s^{\hat m}),L_{\bar A}\wedge L_B\right>\,&= \left<d(\hat
\theta/\hat s^{\hat m}),f_*L_{\bar A}\wedge f_*L_B\right>\\ \,&=
\overline{\gamma^C_A}\gamma^D_B\left<d(\hat \theta/\hat s^{\hat
m}),\hat L_{\bar C}\wedge \hat LD\right>\\ \,&=
-\overline{\gamma^C_A}\gamma^D_B (\hat h^0_{\bar C D}\circ f),
\end{aligned}\end{equation} where the last identity follows from
(\ref{id-2}) and Proposition \ref{m-hAB}. On the other hand, by
(\ref{pull-1}), we have
\begin{equation}\begin{aligned}\Label{calc-12}
\left<df^*(\hat \theta/\hat s^{\hat m}),L_{\bar A}\wedge
L_B\right>\,&=\left<d(\xi\theta/ s^{m}),L_{\bar A}\wedge
L_B\right>\\\ \,&= -\xi h^0_{\bar AB},
\end{aligned}\end{equation}
where again the last identity follows from (\ref{id-2}) and
Proposition \ref{m-hAB}. Thus, we have the identity
\begin{equation}\Label{basic-1}\xi h^0_{\bar A B}=\gamma^D_B\overline
{\gamma^C_A}\hat h^0_{\bar C D}.\end{equation} Here, and in what
follows, we abuse the notation in the following way. For a
function $\hat c$ defined on $\hat M$, we use the notation $\hat
c$ to denote both the function $\hat c\circ f$ on $M$ and the
function $\hat c$ on $\hat M$. It should be clear from the context
which of the two functions is meant. For instance, in
(\ref{basic-1}), we must have $\hat h^0_{\bar C D}=\hat h^0_{\bar
C D}\circ f$. 

By repeating the procedure above to the equation
(\ref{id-1}) with $\hat\omega=\hat\theta/\hat s^{\hat m}$,
$X=L_{\bar A}$, and $Y=S$, we obtain
\begin{equation}\Label{basic-2} L_{\bar A}\xi+\xi h^0_{\bar A}=\xi
\overline{\gamma^C_A}\hat h^0_{\bar C}+\overline{\gamma^C_A}
\eta^D \hat h^0_{\bar C D},\end{equation} where $h^0_{\bar A}$ is
the real-analytic function defined by (\ref{fact-3}). Next,
applying (\ref{id-1}) with $\hat\omega=\hat\theta^E$, $X=L_{\bar
A}$, and $Y=L_B$, we obtain
\begin{equation}\Label{basic-3} L_{\bar A}\gamma^E_B+\eta^Eh^0_
{\bar AB}=0,
\end{equation} by also using the facts that $[L_{\bar A},L_B]$ and
$[\hat L_{\bar C},\hat L_D]$ are multiples of $T$ and $\hat T$
respectively. Applying (\ref{id-1}) with
$\hat\omega=\hat\theta^E$, $X=L_{\bar A}$, and $Y=S$, we obtain
\begin{equation}\Label{basic-4} L_{\bar A}\eta^E+\eta^Eh^0_{\bar
A}=0.
\end{equation} To obtain (\ref{basic-4}), we have used the facts
that commutators of CR vector fields are CR vector fields, and
that $[L_{\bar A},S]$ and $[\hat L_{\bar C},\hat S]$ are multiples
of $T$ and $\hat T$ respectively. Finally, we apply (\ref{id-1})
with $\hat\omega=\hat\theta^E$, $X=S$, and $Y=L_A$ and obtain
\begin{equation}\Label{basic-5}
S\gamma^E_A-L_A\eta^E-\eta^E\overline {h^0_{\bar A}}=0.
\end{equation}
Before proceeding, we observe the following important consequence
of (\ref{basic-1}).

\begin{Pro}\Label{xi} If $M\subset\bC^{n+1}$ is of $m$-infinite
type $2$ at $0\in M$, then the function $\xi$ defined in
$(\ref{push-1})$ is smooth near $0$.
\end{Pro}

\begin{proof} The conclusion follows immediately from
(\ref{basic-1}) and Proposition \ref{m-hAB}.
\end{proof}

Equations (\ref{basic-1})--(\ref{basic-5}) are completely
analogous to (2.9)--(2.13) in \cite{E-UCP}. By following the
arguments in that paper (essentially word for word), repeatedly
applying the vector fields $L_{\bar A}$ to (\ref{basic-1}) and
(\ref{basic-2}), we obtain the following reflection identities
which are analogous to those in \cite{E-UCP}, Theorem 2.4.

\begin{Thm}\Label{reflection} If $\hat M$ is $\hat m$-infinite
$\ell$-nondegenerate at $ 0\in \hat M$, then the following
identities hold for any indices $D,E\in \{1,\ldots, n\}$,
\begin{equation}\Label{ref-id}
\begin{aligned} \gamma^D_E= &\, r_{E}^D\big(\overline{L^J\gamma^C_A},
\overline{L^{ I}\xi};f\big),
\\
\eta^D=&\, s^D\big(\overline{L^J\gamma^C_A}, \overline{L^{
I}\xi};f\big)
\end{aligned}
\end{equation} where \begin{equation}\begin{aligned}\Label{ref-fcn}
&r_{E}^D\big( \overline{L^J\gamma^C_A}, \overline{L^{
I}\xi};q\big)(p),\quad s^D\big(\overline{L^J\gamma^C_A},
\overline{L^{ I}\xi};q\big)(p)
\end{aligned}\end{equation} are real-analytic functions which are
rational in $\overline{L^J\gamma^C_A}$ and polynomial in
$\overline{L^{ I}\xi}$, the indices $A$, $C$ run over the set
$\{1,\ldots, n\}$, and $J$, $I$ over all multi-indices with
$|J|\leq \ell-1$ and $|I|\leq \ell$; here, $(p,q)\in M\times \hat
M$. Moreover, the functions in $(\ref{ref-fcn})$ depend only on
$M$ and $\hat M$ (and not on the mapping $f$).
\end{Thm}
To complete the proof of Theorem \ref{thm-maintech}, we shall use
the following result whose proof follows from that of \cite{E-UCP},
Proposition 3.18.

\begin{Pro}\Label{commutator} If $M\subset\bC^{n+1}$ is of $m$-infinite
type $2$ then, for any multi-index $J$, integer $k\geq 1$, and
index $F\in\{1,\ldots, n\}$ there exist real-analytic functions
$b^{E_1\ldots E_j}_q$  such that \begin{equation}
\Label{commutator-id}
%\multline
\sum_{j=1}^{|J|+k}\sum_{q=0}^{k} b_q^{E_1\ldots E_j}
\underbrace{[\ldots [L_{E_1}\ldots L_{E_j},L_{\bar 1}],L_{\bar
1}]\ldots,L_{\bar 1}]}_{\text{\rm length $q$ }}=L^{J}S^k,
%\endmultline
\end{equation} where standard multi-index notation is used
and the length of the commutator $[\ldots
[X,Y_{1}],Y_{2}]\ldots,Y_{q}]$ is $q$.
\end{Pro}

The arguments in \cite{E-UCP} (indeed, with the simplification
described in the remark following the proof of Theorem 2 due to
the assumption that $M$ is of $m$-infinite type $2$ at $0$) with
$T$ replaced by $S$ now shows that for any multi-indices $R$ and
$Q$, any nonnegative integer $k$, and any indices $D,F\in
\{1,\ldots n\}$, there are smooth functions, which are rational in
their arguments preceeding the ``;'', such that
\begin{equation}\Label{system-1}\begin{aligned} L^RS^kL^{\bar Q}
\gamma^D_F =&\, r^{R\bar Q k}
\big(L^IS^j\gamma^C_A,L^IS^j\eta^C,L^IS^j\xi;f\big),
\\
L^RS^kL^{\bar Q} \eta^D_F =&\, s^{R\bar Q k}
\big(L^IS^j\gamma^C_A,L^IS^j\eta^C,L^IS^j\xi;f\big)\\
L^RS^kL^{\bar Q} \xi=&\, t^{R\bar Q k}
\big(L^IS^j\gamma^C_A,L^IS^j\eta^C,L^IS^j\xi,\overline{L^IS^j\gamma^C_A},
\overline{L^IS^j\eta^C}, \overline{L^IS^j\xi} ;f\big)
\end{aligned}\end{equation} where
$|I|+j\leq 2\ell$.
 The conclusion of Theorem \ref{thm-maintech}
follows by writing (\ref{system-1}), for all $R$, $Q$, $k$ such
that $$|R|+|Q|+k=2\ell+1$$ in the coordinate system $(x,s)$, where
$x=(\re z_1,\im z_1,\ldots, \re z_n,\im z_n)$, for $M$ near $0$
and any coordinate system $\hat y$ for $\hat M$ near $0\in \hat
M$, and observing that the same system of differential equations
holds for any CR mapping $f$ sending a neighborhood of $0$ in $M$
into $\hat M$ with $f(0)$ sufficiently close to $0$. This
completes the proof of Theorem \ref{thm-maintech}.\qed

\section{Proofs of Theorems \ref{Main-1} and
\ref{Main-2}}\Label{pf-Main12}

\begin{proof}[Proof of Theorem $\ref{Main-2}$] By the
Hanges--Treves propagation theorem (see \cite{HT}), it suffices to
show that $f$ extends holomorphically to a full neighborhood of
some particular point $p_1\in E\subset M$. We shall first choose a
point $p_0\in E$
so that Theorem \ref{thm-maintech} is applicable. First, since $M$
is of $1$-infinite type along $E$, $M$ is in fact of $1$-infinite type $2$
outside a proper real-analytic variety of $E$ in view of
Proposition \ref{m-hAB}. Also, since $f|_E\colon E\to \bC^{n+1}$
is finite and $f(E)\subset M'$, $E':=f(E)$ is a
connected complex hypersurface  contained in $M'$. By assumption,
$M'$ is $m'$-essential at some point of $E'$, for some integer
$m'\geq 1$. (It is not difficult to see that, in fact, $m'$ has to
be one.) It follows, exactly as in the finite type case (see e.g.\
\cite{BER}, Chapter IX; cf. also the arguments in \cite{HMM}) that
$M'$ is $m'$-infinite $n$-nondegenerate outside a proper
real-analytic subvariety of $E'$. Since $f|_E$ is finite, we can
find $p_0\in E$ such that $f|_E\colon E\to E'$ is a local
biholomorphism, $M'$ is $m'$-infinite $n$-nondegenerate at
$p_0':=f(p_0)$, and, in addition, $M$ is of $1$-infinite type $2$
at $p_0$. Hence, we may apply Theorem \ref{thm-maintech} to the
mapping $f^0=f$ at $p_0$ with $m=1$ and $\ell=n$.

Let us choose local coordinates $(x,s)\in \bR^{2n}\times \bR$,
vanishing at $p_0$, on $M$ and $y\in \bR^{2n+1}$, vanishing at
$p_0'$, on $M'$ as described in Theorem \ref{thm-maintech}. We
shall denote the components $y_i\circ f$ by $u_i$ in order not to
confuse them with the components of the mapping into $\bC^{n+1}$.
We shall write, for each $i\in\{1,\ldots,2n\}$, each multi-index
$\alpha\in \mathbb Z_+^{2n}$ and each non-negative integer $p$,
$$u_i^{\alpha, p}:=(s\partial_s)^p\partial_x^{\alpha}u_i$$ and also
$U$ for the vector $(u_i^{\beta, q})$, where $|\beta|+q\leq
k:=2n+1$ and $i\in\{1,\ldots,2n\}$. Hence, by Theorem
\ref{thm-maintech}, we have
\begin{equation}\Label{pfaffian-1}
(s\partial_s) u_{i}^{\alpha, p} =r^{\alpha,p}_i(U),\end{equation} for
each $i\in\{1,\ldots,2n\}$ and each multi-index $\alpha$ and
non-negative integer $p$ such that $|\alpha|+p=k$. If we add the
contact equations \begin{equation}\Label{pfaffian-2} (s\partial_s)
u_{i}^{\alpha, p}=u_{i}^{\alpha, p+1},\end{equation} for $|\alpha|+p<k$,
then we obtain the system \begin{equation}\Label{pfaffian-main}
(s\partial_s) U = R(U),\end{equation} where $R(U)(x,s)$ is a
real-analytic vector valued function of $U$ and $(x,s)$,  depending
only on $M$ near $p_0=(0,0)$ and $M'$ near $p_0'=0$ (and the
possibly the value of $U(0,0)$).

Let us fix $x\in \bR^{2n}$ near $0$ and consider
(\ref{pfaffian-main}) as a system of ordinary differential
equations for $U(x,\cdot)$. This system has a singularity of
so-called Briot--Bouquet type at $s=0$, and its properties are well
understood (see e.g. \cite{Sachdev}, Chapters 3.6--3.7 and 8.8;
or \cite{Hille} and further references in these books). We
shall use the following result, which is undoubtedly known. However, 
the author has not found a satisfactory reference for it and, hence,
we shall provide a proof.

\begin{Thm}\Label{reg-thm} Let $y(t)$, with $y=(y_1,\ldots, y_N)$ and
  $t\in \bR$, be a
  $C^\infty$-smooth $\bC^N$-valued function near $t=0$ such that
  $y(0)=0$ and \begin{equation}\Label{BB-system} t\, \frac{dy_j}{dt}=f_j(t,y),\quad
  j=1,\ldots, N,\end{equation} where the $f_j(t,y)$ are analytic functions
  near $(t,y)=(0,0)$. Then, $y(t)$ is real-analytic near $0$.
\end{Thm}

\begin{proof} A
  classical result due to Malmquist (\cite{Malm}; cf. also
  \cite{Hille}, Chapter 11.1) states that the general solution
$y(t)=y_1(t)$ of (\ref{BB-system}), with $N=1$ and
$\lim_{t\to 0} y(t)=0$, is given for $t>0$ by a convergent
series of the form
\begin{equation}\Label{y-series}
y(t)=\sum_{j=0}^\infty \sum_{r=0}^{r_j}c_{ljr} (\ln t)^r
t^{\nu_j},
\end{equation} where the $r_j$ are integers, and $0<\nu_0<\nu_1<\ldots<
\nu_j\nearrow\infty$. A similar convergent series represents
$y(t)$ for $t<0$. If we assume that the solution $y(t)$ is
$C^\infty$-smooth near $t=0$, then no fractional powers or logarithmic
terms can appear and the series for $t<0$ must match that for
$t>0$. We conclude that $y(t)$ is given by a convergent power
series in $t$ and, hence, $y(t)$ is real-analytic near $t=0$. This
completes the proof of Theorem \ref{reg-thm} in the special case
$N=1$. (A similar result for arbitrary $N$ is implicit in the
literature,  but the author has been unable to find an explicit
reference for such a result.)

To treat the case of a general $N$, we shall use a result by Dulac
(building on an idea of Poincare) reducing the system
(\ref{BB-system}) to a special form. After an invertible linear
transformation of the $y_j$ if necessary, we may assume that
$f=(f_1,\ldots, f_N)$ is of the form
\begin{equation}\Label{f-form} f(t,y)=pt+Ay+O(2),\end{equation}
where $p\in\bC^N$, $A$ is an $N\times N$-matrix in Jordan normal
form, and $O(2)$ as usual denotes terms of at least order two in
all the variables. We shall denote the eigenvalues of $A$,
repeated with multiplicity, by $\lambda_1,\ldots, \lambda_N$. The
solution curve $t\to (t,y(t))\in \bC\times\bC^N$ satisfies the differential
system
\begin{equation}\Label{diff-system}
  \frac{dt}{t}=\frac{dy_1}{f_1(t,y)}=\ldots=
  \frac{dy_N}{f_N(t,y)}.
  \end{equation}
  The characteristic roots of
  this system are $1,\lambda_1,\ldots, \lambda_N$.
  By a theorem of
  Dulac (\cite{Dulac}) there are an integer $p$ (determined by the
  location of the roots $1,\lambda_1,\ldots, \lambda_N$ in the
  complex plane), with $0\leq p\leq N$,
  and analytic functions $\Phi_j(t,z_1,\ldots, z_p)$, $j=1,\ldots, p$,
  $\Psi_i(t,z_1,\ldots, z_p)$, $i=1,\ldots , N-p$, with
  $\Phi_j(t,z)=z_j+O(2)$ and $\Psi_i(t,z)=O(2)$ such that
  \begin{equation}\Label{y-z} \begin{aligned} y_j\, &=\Phi_j(t,z),
  \quad j=1,\ldots, p,\\
  y_{p+i}\, &=\Psi_i(t,z),\quad i=1,\ldots,
  N-p\end{aligned}\end{equation} and the system (\ref{diff-system}),
  provided that $p\geq 1$, 
  pulled back to the $(t,z)$-space, becomes
  \begin{equation}\Label{z-system}
  \frac{dt}{t}=\frac{dz_1}{g_1(t,z)}=\ldots=
  \frac{dz_p}{g_p(t,y)},
  \end{equation} where \begin{equation} \Label{g-form}
  g_j(t,z)=\lambda_j
  z_j+P_j(t,z_1,\ldots, z_{j-1}),\quad  j=1,\ldots,
  p,\end{equation} and $P_j$ are polynomials which do not involve
  the variables $z_j,\ldots, z_p$.  (If $p=0$ (which
  corresponds to all eigenvalues $\lambda_1,\ldots, \lambda_N$
  being located on the negative real line $\bR_-\cup\{0\}$), then there
  are no variables $z$ in
  (\ref{y-z}). This already implies that the $y_i(t)$, $i=1,\ldots,
  N$ are real-analytic. Thus, in what follows we may assume that
  $p\geq 1$.)  Observe that the $p$ first
  equations in (\ref{y-z}) can be solved for
  $z_j=\Theta_j(t,y_1,\ldots, y_p)$, $j=1,\ldots,p$. Hence,
  by substituting the smooth functions $y=y(t)$ in this identity,
  we conclude that
  $z=z(t)$ is a $C^\infty$-smooth $\bC^p$-valued function that
  satisfies the system \begin{equation}\Label{BB-z-system} t\,
  \frac{dz_j}{dt}=g_j(t,y),\quad
  j=1,\ldots, p,\end{equation} where $g_j$ are as in
  (\ref{g-form}). By applying the theorem of Malmquist cited above
  to the equation for $z_1$ (which is a single equation of
  Briot-Bouquet type for $z_1(t)$), we conclude that $z_1(t)$
  is real-analytic near $t=0$. By substituting the real-analytic
  function $z_1(t)$ into the equation for $z_2(t)$, using the
  "triangular" form of the equations (\ref{BB-z-system}), and again
  applying Malmquist's theorem, we conclude that $z_2(t)$ is
  real-analytic. By repeating this procedure inductively, we conclude
  that all
  the $z_1(t),\ldots,z_p(t)$ are real-analytic near $0$. The
  real-analyticity of $y(t)$ now follows from (\ref{y-z}). This
  completes the proof of Theorem \ref{reg-thm}.
\end{proof}
By applying Theorem \ref{reg-thm} for fixed $x$ (possibly after
subtracting the value $U(x,0)$ from $U(x,s)$), we conclude that
$U(x,\cdot)$ is given by a convergent series
\begin{equation}\Label{U-series}
U_l(x,s)=\sum_{j\geq 0} a_{lj}(x) s^{j}.
\end{equation}
Moreover, the smoothness of $U$ also implies that the coefficients
$a_{lj}(x)$ are smooth functions of $x$. A simple Baire category
argument shows that there are $x^1\in \bR^{2n}$ near $0$,
$\epsilon>0$, and $\delta>0$ such that the series (\ref{U-series})
converge uniformly in $x$, with $|x-x^1|\leq \epsilon$, for $|s|\leq
\delta$. It is well known (see e.g.\ \cite{BER}, Proposition
1.7.5) that this implies that the CR mapping $f$ extends
holomorphically to a full neighborhood of $p_1=(x^1,0)\in M$ in $\bC^{n+1}$.
This completes the proof of Theorem \ref{Main-2}, in view of the
remark at the beginning of the proof.
\end{proof}

\begin{proof}[Proof of Theorem $\ref{Main-1}$] The proof of Theorem
\ref{Main-1} will follow from Theorem \ref{Main-2} if we can show
that $M'$ is necessarily weakly essential at some point of
$E':=f(E)$. Exactly as in the finite type case, one can easily
show that in $\bC^2$ being $m'$-essential, for some integer $m'$,
at some point is equivalent to being of $m'$-infinite type $r'$
for some integer $r'$ (cf.\ e.g.\ \cite{EH}). Also, as mentioned
in the proof of Theorem \ref{Main-2} above, a real-analytic
hypersurface $M'$ is of $m'$-infinite type $2$ on a dense set of
$E'\subset M'$ if it is of $m'$-infinite type along $E'$. Hence,
to prove Theorem \ref{Main-1} it suffices to show that $M'$ is of
$m'$-infinite type along $E'$ for some integer $m'$; i.e. we must
show that $M'$ is not Levi flat ($m'=\infty$). But this follows
easily from the fact that $M$ is not Levi flat ($m=1$) and condition (i)
(see e.g.\ \cite{E-JGA}, Theorem 2.2). This completes the proof of
Theorem \ref{Main-1}.
\end{proof}

\end{document}